\documentclass[oneside,english]{amsart}
\usepackage[T1]{fontenc}
\usepackage[utf8]{inputenc}
\usepackage{amsthm}
\usepackage{amstext}
\usepackage{amssymb}

\makeatletter
\numberwithin{equation}{section}
\numberwithin{figure}{section}
\theoremstyle{plain}
\newtheorem{thm}{\protect\theoremname}
  \theoremstyle{plain}
  
  \theoremstyle{remark}
  
  \theoremstyle{definition}
  
  \theoremstyle{plain}
  
    \theoremstyle{definition}

\usepackage{babel}
  \providecommand{\definitionname}{Definition}
  \providecommand{\lemmaname}{Lemma}
  \providecommand{\propositionname}{Proposition}
  \providecommand{\remarkname}{Remark}
\providecommand{\theoremname}{Theorem}
\providecommand{\examplename}{Example}

\DeclareMathOperator{\sing}{Sing}

\DeclareMathOperator{\Div}{div}

\author{Rudy Rosas}

\begin{document}

\title{Transversely product singularities of foliations in projective spaces }

\maketitle
\begin{abstract}We prove that a transversely product component of the singular set of  a
holomorphic foliation on $\mathbb P^n$ is necessarily a Kupka component.
\end{abstract}
\section{Introduction}
Let $U$ be an open set of a complex manifold $M$ and let $k\in\mathbb N$. Let $\eta$ be a  holomorphic $k$-form  on $U$ and let
$\sing \eta\colon =\{p\in U\colon \eta(p)=0\}$ denote the singular set of $\eta$.  We say that $\eta$ is integrable if each point $p\in U\backslash \sing \eta$  has a neighborhood $V$ supporting holomorphic  1-forms 
$\eta_1,\dots, \eta_k$ with $\eta|_V=\eta_1\wedge\dots\wedge\eta_k$, such that $d\eta_j\wedge \eta=0$ for each $j=1,\dots, k$. In this case the distribution 
$$\mathcal D_\eta\colon\;\; \mathcal D_\eta (p)=\{v\in T_pM\colon i_v\eta (p)=0\}, \quad 
p\in   U\backslash \sing \eta$$ defines a holomorphic foliation 
of codimension $k$ on  $U\backslash \sing \eta$.
A singular holomorphic foliation $\mathcal F$  of codimension $k$
on  $M$ can be defined by an open covering
$(U_j)_{j\in J}$ of $M$ and a  collection of integrable $k$-forms $\eta_j\in\Omega^k(U_j)$ 
 such that
$\eta_i=g_{ij}\eta_j$ for some $g_{ij}\in\mathcal{O}^*(U_i\cap U_j)$ whenever $U_i\cap U_j\neq \emptyset$. The singular set $\sing \mathcal F$ is the proper analytic subset of $M$ given by the union of the sets $\sing\eta_j$. From now on we only consider foliations $\mathcal F$ such that $\sing\mathcal F$ has no component of codimension one.

 Given a singular holomorphic foliation $\mathcal F$ of codimension $k$ on $M$ as above, the Kupka singular set of $\mathcal F$, denoted by $K(\mathcal F)$,  is the union of the sets 
 $$K(\eta_j)=\{p\in U_j\colon \eta_j(p)=0, d\eta_j (p)\neq 0\}.$$
This set does not depend on the collection $(\eta_j)$ of $k$-forms used to define $\mathcal F$.
It is well known that, given $p\in K(\mathcal F)$,   the germ of $\mathcal F$ at $p$
 is holomorphically equivalent to the product of a one-dimensional foliation with an isolated singularity by a regular foliation of dimension $(\dim \mathcal F -1)$. 
 More precisely, if $\dim M=k+m+1$, there exist a holomorphic vector field
$X=X_1{\partial_{x_1}}+ \dots +X_{k+1}\partial_{x_{k+1}}$
on $\mathbb D^{k+1}$ with a unique singularity at the origin, a neighborhood $V$ of $p$ in $M$ and a biholomorphism 
$\psi \colon V \to \mathbb D^{k+1}\times  \mathbb D^{m}$, $\psi(p)=0$,  which conjugates $\mathcal F$
with the foliation $\mathcal F_X$ of $\mathbb D^{k+1}\times  \mathbb D^{m}$ generated by the commuting
vector fields $X, \partial_{y_1},\dots, \partial_{y_m},$ where $y=(y_1,\dots,y_m)$ are the coordinates in $\mathbb D^m$. If
$\mu=dx_1\wedge\dots \wedge dx_{k+1}$, the foliation $\mathcal F_X$  is
also defined by  the $k$-form $\omega= i_X\mu$ and the Kupka condition $d\omega (0)\neq 0$ is equivalent to the inequality $\Div X(0)\neq 0$.  

Following \cite{LN2021},  we say that $\mathcal F$ is a transversely product at 
  $p\in\sing\mathcal F$ if as above there exist a holomorphic vector field $X$ and a biholomorphism $\psi \colon V \to \mathbb D^{k+1}\times  \mathbb D^{m}$ conjugating $\mathcal F$ with $\mathcal F_X$, except that it is not assumed that $\Div X(0)\neq 0$. We say that $\Gamma$ is a 
  local transversely
  product component of $\sing \mathcal F$ if $\Gamma$ is a compact irreducible component of
  $\sing \mathcal F$ and $\mathcal F$ is a transversely product at each $p\in \Gamma$. In particular, 
   if $\Gamma\subset K(\mathcal F)$ we say that 
  $\Gamma$ is a Kupka component --- for more information about 
  Kupka singularities and Kupka components we refer the reader to \cite{K,CLN, B, CA99, CA09, CA16, CCF}.  If $\Gamma$ is a transversely product component of $\sing \mathcal F$, we can cover 
  $\Gamma$ by finitely many normal coordinates like $\psi$, with the same vector field $X$:
  that is, there exist
   a holomorphic vector field
$X$
on $\mathbb D^{k+1}$ with a unique singularity at the origin and  a covering of
 $\Gamma$  by open sets 
$(V_{\alpha})_{\alpha\in A}$ such that
each $V_\alpha$  supports a biholomorphism 
$\psi_\alpha\colon V_\alpha \to \mathbb D^{k+1}\times  \mathbb D^{m}$ that 
maps $\Gamma\cap V_\alpha$ onto $\{0\}\times\mathbb D^{m}$ and conjugates $\mathcal F$
with the foliation $\mathcal F_X$. The sets $(V_\alpha)$ can be chosen arbitrarily
close to $\Gamma$.

  In 
  \cite{LN2021}, the author proves that a  local transversely
  product component of a codimension one foliation on $\mathbb P^n$ is necessarily a Kupka component. The goal of the present paper is to generalize this theorem to foliations of any codimension. 
  
\begin{thm}\label{teopro}
Let  $\mathcal F$ a holomorphic foliation of dimension $\ge 2$ and codimension $\ge 1$ on 
$\mathbb P^n $.  Let  $\Gamma$ be a transversely product component of $\sing \mathcal F$.
Then  $\Gamma$ is a Kupka component.
\end{thm}

This theorem is a corollary of the following result.
\begin{thm}\label{teotubo}
Let  $\mathcal F$ a holomorphic foliation of dimension $\ge 2$ and codimension $k\ge 1$ on a complex manifold $M$. Suppose that $\mathcal F$ is defined by an open covering
$(U_j)_{j\in J}$ of $M$ and a  collection of $k$-forms $\eta_j\in\Omega^k(U_j)$. Let $L$ be the line bundle defined
by the cocycle $(g_{ij})$ such that 
$\eta_i=g_{ij}\eta_j$, $g_{ij}\in\mathcal{O}^*(U_i\cap U_j)$.   Let  $\Gamma$ be a transversely product component of $\sing \mathcal F$ that is not a Kupka
component. Then, if $V$ is a  tubular neighborhood of $\Gamma$, we have that  $c_1(L|_V)=0$.
\end{thm}

\section {Proof of the results}
\noindent{\emph{Proof of Theorem \ref{teotubo}.} }
Let $V$ be  a tubular neighborhood of $\Gamma$. Then the map $$\Theta\in H^2_{\textrm{dR}}(V)\mapsto \Theta|_\Gamma \in H^2_{\textrm{dR}}(\Gamma)$$ is an isomorphism and so it suffices
to prove that $c_1(L|_\Gamma)=0$.   
Let  $\dim M=k+m+1$. As explained in the introduction, there exist a holomorphic vector field
$X$
on $\mathbb D^{k+1}$ with a unique singularity at the origin  and  a covering of
 $\Gamma$  by open sets 
$(V_{\alpha})_{\alpha\in A}$ such that
each $V_\alpha$ is contained in $V$ and supports a biholomorphism 
$\psi_\alpha\colon V_\alpha \to \mathbb D^{k+1}\times  \mathbb D^{m}$ that 
maps $\Gamma\cap V_\alpha$ onto $\{0\}\times\mathbb D^{m}$ and conjugates $\mathcal F$
with the foliation $\mathcal F_X$  generated by the commuting
vector fields $X, \partial_{y_1},\dots, \partial_{y_m}$. Notice that  $\Div(X)(0)=0$, because $\Gamma$ is not a Kupka component.
  Since $\mathcal F_X$ is defined  by the  $k$-form $\omega= i_X\mu$, where
$\mu=dx_1\wedge\dots\wedge dx_{k+1}$, 
 we have
 that
  $\mathcal F|_{V_\alpha}$ is defined by
the $k$-form $\psi^*_\alpha(\omega)$. If $V_\alpha\cap V_\beta\neq\emptyset$, there exists
 $\theta_{\alpha\beta}\in\mathcal{O}^*(V_\alpha\cap V_\beta)$ such that 
 \begin{align}\label{cociclo} \psi^*_\alpha(\omega)=\theta_{\alpha\beta}\psi^*_\beta(\omega).
 \end{align}
  We can assume 
 that the $k$-forms $\psi^*_\alpha(\omega)$ belong to the family of $k$-forms $(\eta_j)_{j\in J}$
 defining $\mathcal F$. Therefore the cocycle $(\theta_{\alpha\beta})$ define the line bundle $L$ restricted to some neighborhood of $\Gamma$. Thus, in order to prove that $c_1(L|_\Gamma)=0$ 
 it is enough to show that each $\theta_{\alpha\beta}|_\Gamma$ is locally constant. Fix some 
 $\alpha, \beta\in A$ such that $V_\alpha\cap V_\beta\neq \emptyset$. If we set
  $\psi=\psi_\alpha\circ\psi_\beta^{-1}$ and $\theta=\theta_{\alpha\beta}\circ\psi_\beta^{-1}$, from \eqref{cociclo}
 we have that 
$\psi^*(\omega)=\theta \omega,
$ which means that $\psi$ preserves the foliation $\mathcal F_X$.
 It suffices to prove that the derivatives
 $\theta_{y_1}(p),\dots,  \theta_{y_m}(p)$ vanish
 if
 $p\in \{0\}\times \mathbb D^m$. 
 Since $\partial_{y_1}$ is tangent to $\mathcal F_X$, then the vector field ${\psi_*({\partial_{y_1}})}$ is tangent to
  $\mathcal F_X$ and so  we can express
 $$\psi_*({\partial_{y_1}})=\lambda X +\lambda_1\partial_{y_1}+\dots + \lambda_m\partial_{y_m},$$
 where $\lambda,\lambda_1,\dots,\lambda_m$ are holomorphic. Then
 \begin{align*}\label{ecua1}\mathcal{L}_{\psi_*({\partial_{y_1}})}\omega=\mathcal{L}_{\lambda X}\omega=
 \lambda\mathcal{L}_{X}\omega+d\lambda\wedge i_X \omega= \lambda\mathcal{L}_{X}\omega
 =\lambda \Div (X)\omega,
 \end{align*} where the last equality follows from the identity $\omega=i_X\mu$. 
Thus, since 
 $$\psi^*\left(\mathcal{L}_{\psi_*({\partial_{y_1}})}\omega\right)=
 \mathcal{L}_{\partial_{y_1}}\psi^*\omega=\mathcal{L}_{\partial_{y_1}}(\theta\omega)
 =\theta_{y_1}\omega,$$ we obtain that
 $$\theta_{y_1}\omega=\psi^*\left(\lambda \Div (X)\omega\right)=
 \lambda (\psi) \Div (X)(\psi) \theta\omega$$
 and therefore $\theta_{y_1}(p)=0$ if 
  $p\in \{0\}\times \mathbb D^m$, because
 $\Div (X)$ vanishes along $\{0\}\times \mathbb D^m$. In the same way we prove that 
 $\theta_{y_2}(p)=\dots =\theta_{y_m}(p)=0$  if 
  $p\in \{0\}\times \mathbb D^m$, which finishes the proof. \qed \\

 \noindent{\emph{Proof of Theorem \ref{teopro}.} }
Suppose that $\Gamma$ is not a Kupka component. Let $L$ be the line bundle associated to $\mathcal F$ as in the statement of Theorem \ref{teotubo}. 
We notice that $c_1(L)\neq 0$, otherwise $\mathcal F$ will be defined by a global $k$-form on 
$\mathbb P^n$, which is impossible. Then, if we take an algebraic curve 
$\mathcal C \subset \Gamma$, we have 
$c_1(L)\cdot\mathcal C\neq 0$.  Therefore, if $\Omega$ is a 2-form on $\mathbb P^n$ in the class $c_1(L)$ and $V$ is a tubular neighborhood of $\Gamma$, 
$$c_1(L|_V)\cdot\mathcal C =\int_{\mathcal C}\Omega|_V=\int_{\mathcal C}\Omega=c_1(L)\cdot
\mathcal C\neq 0,$$ which contradicts Theorem \ref{teotubo}. \qed

\end{document}